\documentclass{amsart}
\usepackage{amssymb}
\usepackage{amsmath,amsthm,amscd}
\usepackage[mathscr]{eucal}

\newcommand{\Slo}{\mathcal{S}}
\newcommand{\pause}{}
\newcommand{\Z}{\mathbb{Z}}
\newcommand{\W}{\mathscr{W}}
\newcommand{\ra}{\rightarrow}
\newcommand{\on}{\operatorname}
\newcommand{\+}{\oplus}
\newcommand{\C}{\mathbb{C}}
\newcommand{\g}{\mathfrak{g}}
\newcommand{\affg}{\widehat{\mathfrak{g}}}
\newcommand{\mf}{\mathfrak}
\newcommand{\mc}{\mathcal}
\newcommand{\Lam}{\Lambda}
\newcommand{\lam}{\lambda}

\newcommand{\isomap}{{\;\stackrel{_\sim}{\to}\;}}

\begin{document}
\theoremstyle{theorem}
\newtheorem{Th}{Theorem}
\newtheorem*{Pro}{Proposition}
\newtheorem{Lem}[Th]{Lemma}
\theoremstyle{definition}

\newtheorem*{Prob}{基本問題}
\newtheorem{Rem}[Th]{Remark}
\newtheorem{Fact}{事実}
\newtheorem{Conj}[Th]{Conjecture}

\newtheorem*{Def}{Definition}
\newtheorem{Co}[Th]{系}
\theoremstyle{remark}
\newtheorem{MainTh}{Main Theorem}

\newtheorem*{Facts}{Fact}

\newtheorem*{FP}{Fundamental Problem}
\newtheorem*{Goal}{Today's Goal}

\newtheorem*{Note}{Note}

\newtheorem*{Open}{Open Problem}
\newtheorem{Claim}{Claim}
\newtheorem*{ClaimNN}{Claim}
\newtheorem{Notation}{Notation}
\newtheorem*{Question}{Question}
\newtheorem*{Ex}{例}

\renewcommand{\baselinestretch}{1.1}

\title{Representation theory of W-algebras
and Higgs branch conjecture}

\author{Tomoyuki Arakawa}
\address{Research Institute for Mathematical Sciences, Kyoto University,
Kyoto 606-8502 JAPAN}
\address{Department of Mathematics,
Massachusetts Institute of Technology,
77 Massachusetts Avenue
Cambridge, MA 02139-4307 USA}
\email{arakawa@kurims.kyoto-u.ac.jp}
\maketitle

\begin{abstract}
We survey a number of  results regarding the representation theory of $W$-algebras
and their connection with the resent development of the four dimensional $N=2$ superconformal field theories
in physics.
\end{abstract}

 
\section{Introduction}
(Affine)  $W$-algebras appeared in 80's in the study of the two-dimensional conformal field theory in physics.
$W$-algebras can
 be regarded as a generalization of infinite-dimensional Lie algebras such as 
affine Kac-Moody algebras and the Virasoro algebra,
although $W$-algebras are not Lie algebras but vertex algebras in general. 
$W$-algebras may be also considered as an affinizaton  of 
finite $W$-algebras (\cite{BoeTji93,Pre02})
which are
a natural quantization of Slodowy slices.
 $W$-algebras play
important roles not only in conformal field theories
but also in integrable systems (e.g.\ \cite{DriSok84,De-KacVal13,BakMil13}),
the geometric Langlands program (e.g.\ \cite{Fre07,AgaFreOko})
and four-dimensional gauge theories (e.g. \cite{AGT,SchVas13,MauOko1211,BraFinNak16}).

In this note we survey the resent development of the representation theory of $W$-algebras.
One of the fundamental problems in $W$-algebras was the Frenkel-Kac-Wakimoto  conjecture \cite{FKW92}
that stated the existence and  construction of rational $W$-algebras, 
which generalizes the integrable representations of affine Kac-Moody algebras and  the minimal series representations of the Virasoro algebra.
The notion of the associated varieties of vertex algebras played a crucial role in the proof
\cite{A2012Dec} of the Frenkel-Kac-Wakimoto
conjecture, and  has revealed  \cite{BeeRas} an unexpected connection 
of vertex algebras with the geometric invariants called the {Higgs branches} in the {four dimensional} $N=2$ superconformal field theoriesin physics.

\subsection*{Acknowledgments}
The author benefited greatly from discussion with Christopher Beem, 
Davide Gaiotto,
Madalena Lemos,
Victor Kac, 
Anne Moreau, Hiraku Nakajima, Takahiro Nishinaka, Wolfger Peelaers, Leonardo Rastelli, Shu-Heng Shao, Yuji Tachikawa, and Dan Xie.
The author 
is partially supported in part by JSPS KAKENHI Grant Numbers 17H01086, 17K18724.

\section{Vertex algebras}
A \emph{vertex algebra} 
\cite{Bor86}
consists of a vector space $V$ with a distinguished vacuum vector $|0\rangle  \in V$ and a vertex operation, which is a linear map $V \otimes V \rightarrow V((z))$, written $a \otimes b \mapsto a(z)b = (\sum_{n \in \Z} a_{(n)} z^{-n-1})b$, such that the following are satisfied:
\begin{itemize}
\item (Unit axioms) $(|0\rangle)(z) = 1_V$ and $a(z)|0\rangle  \in a + zV[[z]]$ for all $a \in V$.

\item (Locality)
$(z-w)^n[a(z),b(w)]=0$ for a sufficiently large $n$
for all $a, b\in V$.
\end{itemize}
The operator $T: a \mapsto a_{(-2)}|0\rangle $ is called the {\em translation operator} and it satisfies $(Ta)( z) =
[T,a(z)]= \partial_z a(z)$. The operators $a_{(n)}$ are called \emph{modes}.

For elements $a,b$ of a vertex algebra $V$
 we have the following
{\em Borcherds identity}
for any
$m,n\in \Z$:
\begin{align}
&[a_{(m)},b_{(n)}]=\sum_{j\geq 0}\begin{pmatrix}m\\j\end{pmatrix}(a_{(j)}b)_{(m+n-j)},
\label{eq:B1}\\
& (a_{(m)}b)_{(n)}=\sum_{j\geq 0}(-1)^j\begin{pmatrix}m\\j\end{pmatrix}
(a_{(m-j)}b_{(n+j)}-(-1)^m b_{(m+n-j)}a_{(j)}).
\end{align}
By regarding the Borcherds identity
as fundamental relations,
 representations of a vertex algebra are naturally defined
(see \cite{Kac98,FreBen04} for the details).

\smallskip

One of the basic examples of vertex algebras
are universal affine vertex algebras.
Let $G$ be a simply connected simple algebraic group,
$\g=\on{Lie}(G)$.
Let 
$\widehat{\g}=\g[t,t^{-1}]\+ \C K$
be
the 
affine Kac-Moody algebra associated with $\g$.
The commutation relations of $\affg$ are given by
\begin{align}
[xt^m, yt^n]=[x,y]t^{m+n}+m\delta_{m+n,0}(x|y)K,\quad [K,\affg]=0 \quad (x,y\in \g,\ m,n\in \Z),
\label{eq:affg-com-rel}
\end{align}
where
$(~|~)$ is the normalized invariant inner product of $\g$,
that is,
$(~|~)=1/2h^{\vee}\times $Killing form
and  $h^{\vee}$ is the dual Coxeter number of $\g$.
For $k\in \C$,
let
\begin{align*}
V^k(\g)=U(\affg)\otimes_{U(\g[t]\+ \C K)}\C_k,
\end{align*}
where 
$\C_k$ is the one-dimensional representation
of $\g[t]\+ \C K$
on which 
$\g[t]$ acts trivially and
$K$ acts as multiplication by $k$.
There is a unique vertex algebra structure
on 
$V^k(\g)$
such that $|0\rangle=1\otimes 1$ is the vacuum vector
and
\begin{align*}
x(z)=\sum_{n\in \Z}(xt^n)z^{-n-1}\quad (x\in \g).
\end{align*}
Here on the left-hand-side
$\g$ is considered as a subspace 
of $V^k(\g)$ by the embedding
$\g\hookrightarrow V^k(\g)$,
$x\mapsto (xt^{-1})|0\rangle$.
$V^k(\g)$
is called the {\em universal affine vertex algebra associated with $\g$ at a level $k$}.
The Borcherds identity \eqref{eq:B1} for $x,y\in \g\subset V^k(\g)$
is identical to the commutation relation \eqref{eq:affg-com-rel}
with $K=k\on{id}$,
and hence, any $V^k(\g)$-module
is a $\affg$-module of level $k$.
Conversely,
any smooth $\affg$-module of level $k$
is naturally
a $V^k(\g)$-module,
and therefore, the category 
$V^k(\g)\on{-Mod}$ of $V^k(\g)$-modules is the same as the category of smooth $\affg$-modules
of level $k$.
Let $L_k(\g)$ be the unique simple graded quotient of $V^k(\g)$,
which is isomorphic to the irreducible highest weight representation $L(k\Lam_0)$ with highest weight 
$k\Lam_0$ as a $\affg$-module.
The vertex algebra $L_k(\g)$ is called the {\em simple affine vertex algebra} associated with $\g$ at level $k$,
and $L_k(\g)\on{-Mod}$ forms a full subcategory of $V^k(\g)\on{-Mod}$,
 the category of smooth $\affg$-modules
of level $k$.

\smallskip

A vertex algebra $V$ is called {\em commutative}
if both sides of \eqref{eq:B1} are zero for all $a,b\in V$,
$m,n\in \Z$.
If this is the case, $V$ can be regarded as a {\em differential algebra} (=a unital commutative algebra with a derivation)
 by the multiplication
$
a.b=a_{(-1)}b
$
and the derivation
$T$.
Conversely,
any differential algebra can be naturally equipped with the structure of a commutative vertex algebra.
Hence,
commutative vertex algebras are the same\footnote{However,
the modules of a commutative vertex algebra are not the same as the modules
as a differential algebra.} as differential algebras (\cite{Bor86}).

Let $X$ be an affine scheme,
$J_{\infty}X$ the
{\em arc space} of $X$ that is defined by the functor of points
$\on{Hom}(\on{Spec}R,J_{\infty}X)=\on{Hom}(\on{Spec}R[[t]], X)$.
The ring $\C[J_{\infty}X]$ is naturally a differential algebra,
and hence is a commutative vertex algebra.
In the case that 
$X$ is a Poisson scheme
 $\C[J_{\infty}X]$ has \cite{Ara12} the structure of 
{\em Poisson vertex algebra},
which is a vertex algebra analogue of Poisson algebra (see \cite{FreBen04,Kaclecture} for the precise definition).

It is known by  Haisheng Li \cite{Li05} that
any vertex algebra $V$ is canonically filtered,
and hence can be regarded\footnote{This filtration is separated if $V$ is non-negatively graded,
which we assume.} as a quantization of 
the associated graded Poisson vertex algebra 
$\on{gr} V=\bigoplus_p F^pV/F^{p+1}V$,
where 
$F^{\bullet}V$ is the canonical filtration of $V$.
By definition,
\begin{align*}
F^pV=\on{span}_\C\{(a_1)_{(-n_1-1)}\dots (a_r)_{(-n_r-1)}|0\rangle\mid a_i\in V,\ n_i\geq 0,\ \sum_i n_i\geq p\}.
\end{align*}
The subspace 
 $$R_V:=V/F^1V=F^0V/F^1V\subset \on{gr}V$$
is called 
 {\em Zhu's $C_2$-algebra of $V$}.
The Poisson vertex algebra structure of $\on{gr} V$ restricts to the Poisson algebra structure
of $R_V$,
which  is given by
\begin{align*}
\bar a.\bar b=\overline{a_{(-1)}b},\quad
\{\bar a,\bar b\}=\overline{a_{(0)}b}.
\end{align*}
The Poisson variety
\begin{align*}
X_V=\on{Specm}(R_V)
\end{align*}
called the {\em associated variety}  of $V$ (\cite{Ara12}).
We have \cite{Li05} the inclusion
\begin{align}
\on{Specm}(\on{gr}V)\subset J_{\infty}X_V.
\label{eq:inclusion-SS}
\end{align}

A vertex algebra $V$ is 
called {\em finitely strongly generated} if $R_V$ is finitely generated.
In this note all vertex algebras are assumed to be finitely strongly generated.
$V$ is 
called {\em lisse} (\cite{BeiFeiMaz}) (or {\em $C_2$-cofinite} (\cite{Zhu96}))  if 
$\dim X_V=0$.
By \eqref{eq:inclusion-SS},
it follows that
$V$ is lisse
if and only if  $ \dim \on{Spec}(\on {gr} V)=0$  (\cite{Ara12}).
Hence lisse vertex algebras can be regarded as an analogue of finite-dimensional algebras.

For instance, 
consider the case $V=V^k(\g)$.
We have
$F^1V^k(\g))=\g[t^{-1}]t^{-2}V^k(\g)$,
and there is  an  isomorphism of Poisson algebras
\begin{align*}
\C[\g^*]\isomap  R_V,\quad x_1\dots x_r\mapsto \overline{(x_1t^{-1})\dots (x_rt^{-1})|0\rangle}\quad(x_i\in \g).
\end{align*}
Hence
\begin{align}
X_{V^k(\g)}\cong \g^*.
\label{eq:ass-var--affine}
\end{align}
Also,
we have the isomorphism 
$\on{Spec}(\on{gr}V^k(\g))\cong J_{\infty}\g^* $. 
By \eqref{eq:ass-var--affine},
we have $X_{L_k(\g)}\subset \g^*$,
which is $G$-invariant and conic.
It is known \cite{DonMas06} that 
\begin{align}
\text{$L_k(\g)$ is lisse }\iff \text{ $L_k(\g)$ is integrable as a $\affg$-module}
\text{ ($\iff$ $k\in \mathbb{Z}_{\geq 0}$)}.
\label{eq:integrability}
\end{align}
Hence
the lisse condition may be regarded as a generalization of the integrability condition to an arbitrary 
vertex algebra.

\smallskip

A vertex algebra is called {\em conformal}
if there exists a vector $\omega$,
called the {\em conformal vector},
such that the corresponding field $\omega(z)=\sum_{n\in \Z}L_nz^{-n-2}$
satisfies the following conditions. 
(1) $[L_m, L_n]=(m-n)L_{m+n}+\frac{m^3-m}{12}\delta_{m+n,0}c\on{id}_V$,
 where $c$ is a constant called the {\em central charge}  of V;
 (2)
$L_0$ acts semisimply on V;
(3) $L_{-1}=T$.
For a conformal vertex algebra
$V$ we set $V_{\Delta}=\{v\in V\mid L_0v=\Delta V\}$,
so that $V=\bigoplus_{\Delta}V_{\Delta}$.
The universal affine vertex algebra $V^k(\g)$ is conformal by the Sugawara construction
provided that $k\ne -h^{\vee}$.
A {\em positive energy representation} $M$
of a conformal vertex algebra $V$ is a $V$-module  $M$
on which $L_0$ acts diagonally and the $L_0$-eivenvalues on $M$ is bounded from below.
An  {\em ordinary  representation} is a positive energy representation such that
each $L_0$-eivenspaces are finite-dimensional.
For a finitely generated ordinary representation $M$,
the normalized character
\begin{align*}
\chi_V(q)=\on{tr}_V (q^{L_0-c/24})
\end{align*}
is well-defined.

For a conformal vertex algebra $V=\bigoplus_{\Delta}V_{\Delta}$,
one defines {\em Zhu's algebra} \cite{FreZhu92} $\on{Zhu}(V)$ of  $V$ by
\begin{align*}
\on{Zhu}(V)=V/V\circ V,\quad V\circ V=\on{span}_{\C}\{a\circ b\mid a,b\in V\},
\end{align*}
where $a\circ b=\sum_{i\geq 0}\begin{pmatrix}
\Delta_a\\
i\end{pmatrix}a_{(i-2)}b$ for $a\in V_{\Delta_a}$.
$\on{Zhu}(V)$ is a unital associative algebra by the multiplication
$a* b=\sum_{i\geq 0}\begin{pmatrix}
\Delta_a\\
i\end{pmatrix}a_{(i-1)}b$.
There is a bijection between the isomorphism classes
$ \on{Irrep}(V)$
of simple positive energy representation of $V$
and that of simple
$\on{Zhu}(V)$-modules (\cite{FreZhu92,Zhu96})\footnote{In fact, there is an isomorphicm 
$\on{Zhu}(V)\cong \mc{U}(V)/\overline{\bigoplus_{n>0}\mc{U}(V)_{-n}\mc{U}(V)_n}$ (\cite{FreZhu92,NT05, He17}),
where $\mc{U}(V)$ is the current algebra of $V$.}.
The grading of $V$ gives a filtration on $\on{Zhu}(V)$
which makes it 
 quasi-commutative,
 and
there is a surjective map
\begin{align}
R_V\twoheadrightarrow \on{gr}\on{Zhu}(V)
\label{eq:surj-RV}
\end{align}
of Poisson algebras.
Hence,
if $V$ is lisse  $\on{Zhu}(V)$ is finite-dimensional,
and so  there are only finitely many irreducible positive energy representations of $V$.
Moreover, the lisse condition implies that
 any simple $V$-module  is a positive energy representation
(\cite{AbeBuhDon04}).

A conformal vertex algebra is called {\em rational} if any positive energy representation of $V$ is completely reducible.
For instance,
the simple affine vertex algebra $L_k(\g)$ is rational if and only if 
$L_k(\g)$ is integrable,
and if this is the case $L_k(\g)\on{-Mod}$ is exactly the category of integrable representations of $\affg$
at level $k$.
A theorem of  Y. Zhu \cite{Zhu96} states that
if $V$ is a rational, lisse, $\Z_{\geq 0}$-graded conformal vertex  algebra
such that $V_0=\C|0\rangle$, then
the character $\chi_M(e^{2\pi i\tau})$ converges 
to a holomorphic functor on the upper half plane
for 
any $M\in \on{Irrep}(V)$.
Moreover,
the space 
spanned by the characters $\chi_M(e^{2\pi i\tau})$, $M\in \on{Irrep}(V)$,
is invariant under the natural action of $SL_2(\Z)$.
This theorem was strengthened  in \cite{DonLinNg15}  to the fact that
$\{\chi_M(e^{2\pi i\tau})\mid M\in \on{Irrep}(V)\}$  forms a vector valued modular function
by showing the congruence property.
Furthermore, it has been shown in \cite{Hua08rigidity}
that
the category of $V$-modules form a modular tensor category.

\section{ $W$-algebras}
$W$-algebras are defined by
the method of the {\em quantized Drinfeld-Sokolov reduction}
that
was discovered by Feigin and Fenkel
\cite{FF90}.
In the most general definition of $W$-algebras  given by 
Kac and Wakimoto \cite{KacRoaWak03},
$W$-algebras are associated with
the pair $(\g,f)$ of 
 a simple Lie algebra $\g$
 and
a nilpotent element $f\in \g$.
The
corresponding $W$-algebra
is a one-parameter family of vertex algebra 
denoted by $\W^k(\g,f)$, $k\in \C$.
By definition,
$$\W^k(\g,f):=H^0_{DS,f}(V^k(\g)),$$
where 
$H^\bullet_{DS,f}(M)$  denotes the BRST cohomology 
of the quantized Drinfeld-Sokolov reduction associated with $(\g,f)$ with  coefficient in  a $V^k(\g)$-module
$M$,
which is defined as follows.
Let $\{e,h,f\}$ be an $\mf{sl}_2$-triple associated with $f$,
$\g_j=\{x\in \g\mid [h,x]=2j\}$, so that
$\g=\bigoplus_{j\in \frac{1}{2}\Z}\g_j$.
Set $\g_{\geq 1}=\bigoplus_{j\geq 1}\g_j$,
$\g_{>0}=\bigoplus_{j\geq 1/2}\g_j$.
Then $\chi:\g_{\geq 1}[t,t^{-1}]\ra \C$,
$xt^n\mapsto \delta_{n,-1}(f|x)$,
defines a character.
Let $F_{\chi}=U(\g_{>0}[t,t^{-1}])\otimes_{U(\g_{>0}[t]+\g_{\geq 1}[t,t^{-1}])}\C_{\chi}$,
where $\C_{\chi}$ is the one-dimensional representation of 
$\g_{>0}[t]+\g_{\geq 1}[t,t^{-1}]$
on which $\g_{\geq 1}[t,t^{-1}]$ acts by  the character $\chi$
and  $\g_{>0}[t]$ acts triviality.
Then, for a $V^k(\g)$-module $M$, 
\begin{align*}
H^\bullet_{DS,f}(M)=H^{\frac{\infty}{2}+\bullet}(\g_{>0}[t,t^{-1}],M\otimes  F_{\chi}),
\end{align*}
where $H^{\frac{\infty}{2}+\bullet}(\g_{>0}[t,t^{-1}],N)$ is the semi-infinite  $\g_{>0}[t,t^{-1}]$-cohomology
\cite{Feu84}
with coefficient in a $\g_{>0}[t,t^{-1}]$-module $N$.
Since it is defined by a BRST cohomology,
$\W^k(\g,f)$ is naturally a vertex algebra,
which 
is called
the
 {\em W-algebra
associated with $(\g,f)$ at level $k$}.
By \cite{FreBen04,KacWak04},
we know that  $H^i_{DS,f}(V^k(\g))=0$ for $i\ne 0$.
If $f=0$ we have by definition $\W^k(\g,f)=V^k(\g)$.
The $W$-algebra $\W^k(\g,f)$  is conformal provided that $k\ne -h^{\vee}$.

Let $\mc{S}_f=f+\g^e\subset \g\cong \g^*$, the {\em Slodowy slice} at $f$,
where $\g^e$ denotes the centralizer of $e$ in $\g$. 
The affine variety $\mc{S}_f$ has a Poisson structure obtained from that of $\g^*$ by Hamiltonian reduction (\cite{GanGin02}).
We have
\begin{align}
{X}_{\W^k(\g,f)}\cong \mc{S}_f,\quad \on{Spec}(\on{gr}\W^k(\g,f))\cong J_{\infty}\mc{S}_f
\label{eq:Iso-with-Slodowy}
\end{align}
 (\cite{De-Kac06,Ara09b}).
 Also,
 we have
 \begin{align*}
\on{Zhu}(\W^k(\g,f))\cong U(\g,f)
\end{align*}
(\cite{Ara07,De-Kac06}),
where  $U(\g,f)$ is the {\em finite $W$-algebra} associated with $(\g,f)$ (\cite{Pre02}).
Therefore, the $W$-algebra
$\W^k(\g,f)$ can be regarded as an affinization of the finite $W$-algebra $U(\g,f)$.
The map
\eqref{eq:surj-RV}
 for $\W^k(\g,f)$ is an isomorphism (\cite{De-Kac06,A2012Dec}),
 which recovers the fact \cite{Pre02,GanGin02} that 
 $U(\g,f)$ is a 
 quantization of the Slodowy slice $\mc{S}_f$.
The definition of
 $\W^k(\g,f)$ naturally extends \cite{KacRoaWak03} to the case that   $\g$ is a (basic classical) Lie superalgebra and 
 $f$ is a nilpotent element  in the even part of $\g$.

We have $\W^k(\g,f)\cong \W^k(\g,f')$ 
if $f$ and $f'$ belong to the same nilpotent orbit of $\g$.
The $W$-algebra associated with a minimal nilpotent element $f_{mim}$ and a principal nilpotent element $f_{prin}$
are called a minimal $W$-algebra and a principal $W$-algebra, respectively.
For $\g=\mf{sl}_2$, these two coincide and are isomorphic to the Virasoro vertex algebra
of central charge $1-6(k-1)^2/(k+2)$ provided that $k\ne -2$.
In \cite{KacRoaWak03}
it was
shown that 
almost every superconformal algebra
appears  as the minimal $W$-algebra  $\W^k(\g,f_{min})$
 for some
Lie superalgebra $\g$,
by describing the generators and the relations (OPEs) of minimal $W$-algebras.
Unfortunately,
the presentation  of $\W^k(\g,f)$ by generators and relations 
is not known for other nilpotent elements
except for some special cases.

Historically,
the principal $W$-algebras were first extensively studied (see \cite{Bou95} and references therein).
In the case that $\g=\mf{sl}_n$,
the non-critical principal $W$-algebras are isomorphic to the 
Fateev-Lukyanov's $W_n$-algebra \cite{FatLyk88}
 (\cite{FF90,FreBen04}).
The critical principal $W$-algebra $\W^{-h^{\vee}}(\g,f_{prin})$ is  isomorphic to the 
Feigin-Frenkel center  $\mf{z}(\affg)$ of $\affg$,
that is 
the 
center of the critical affine vertex algebra $V^{-h^{\vee}}(\g)$ (\cite{FeiFre92}).
This should be viewed as a chiralization of Kostant's theorem \cite{Kos78}
which states that $U(\g,f_{prin})\cong \mc{Z}(\g)$,
where $\mc{Z}(\g)$ is the center of $U(\g)$.
For a general $f$,
we have  \cite{Pre07} the isomorphism
$ \mc{Z}(\g)\cong Z(U(\g,f))$ for finite $W$-algebras
which upgrades \cite{A11,AMcenter} to the fact that
\begin{align}
\mf{z}(\affg)\cong Z(\W^{-h^{\vee}}(\g,f)),
\label{eq:center}
\end{align}
where $Z(U(\g,f))$ and $Z(\W^{-h^{\vee}}(\g,f))$ are  the center of 
$U(\g,f)$ and $\W^{-h^{\vee}}(\g,f)$, respectively.
The isomorphism \eqref{eq:center}
 has an application to {\em Vinberg's Problem} for the centralizer $\g^e$ of $e$ in $\g$ \cite{AP}.

\section{Representation theory of $W$-algebras}
The definition of $\W^k(\g,f)$ by the quantized Drinfeld-Sokolov reduction 
gives rise to  a functor
\begin{align*}
V^k(\g)\on{-Mod}&\ra \W^k(\g,f)\on{-Mod}\\
M&\mapsto H^0_{DS,f}(M).
\end{align*}
Let 
$\mc{O}_k$
be the category $\mc{O}$ of $\affg$ at level $k$.
Then $\mc{O}_k$ is naturally considered as a full subcategory of 
$V^k(\g)\on{-Mod}$.
For a weight $\lam$ of $\affg$ of level $k$,
let $L(\lam)$ be the  irreducible
highest weight representations of $\affg$
with highest weight $\lam$.
\begin{Th}[\cite{Ara05}]\label{Th:minimal}
Let $f_{min}\in \mathbb{O}_{min}$
and  let $k$ be an arbitrary complex number.
\begin{enumerate}
\item  We have
$H^{i}_{DS,f_{min}}(M)=0$
for any
 $M\in \mc{O}_k$ and $i\in \Z\backslash \{0\}$.
Therefore,
the functor
$$\mc{O}_k\ra \W^k(\g,f_{min})\on{-Mod},\quad
M\mapsto H^0_{DS,f_{min}}(M)$$
is exact.

\item 
For a weight $\lam$ of $\affg$ of level $k$, 
$H^{0}_{DS,f_{min}}(L(\lam))$
is zero or isomorphic to an
irreducible highest weight representation of
$\W^k(\g,f_{min})$.
Moreover, any 
irreducible highest weight representation
the minimal $W$-algebra $\W^k(\g,f_{min})$ arises in this way.
 \end{enumerate}
\end{Th}

By Theorem \ref{Th:minimal}
and the
Euler-Poincar\'{e} principal,
the character
$\on{ch}H^0_{DS,f_{\theta}}(L(\lam))$ 
is expressed in terms of  $\on{ch} L(\lam)$.
Since $\on{ch} L(\lam)$ is known \cite{KasTan00}
for all non-critical weight $\lam$,
 Theorem \ref{Th:minimal}
determines the character of all non-critical 
irreducible highest weight representation of
$\W^k(\g,f_{min})$.
In the case that $k$ is critical
the character of 
irreducible highest weight representation of
$\W^k(\g,f_{min})$ is determined by the Lusztig-Feigin-Frenkel conjecture (\cite{Lus91,AraFie08,FreGai09}).

\pause
\begin{Rem}
Theorem \ref{Th:minimal}
holds in the case that $\g$ is a basic classical Lie superalgebra as well.
In particular
one obtains the character of irreducible highest weight representations of
superconformal algebras  that appear as $\W^k(\g,f_{min})$
once the character of irreducible highest weight representations of
 $\affg$ is  given.
\end{Rem}

Let $\on{KL}_k$ be the full subcategory of $\mc{O}_k$
consisting of objects on which $\g[t]$ acts locally finitely.
Altough
the functor
\begin{align}
\mc{O}_k\ra \W^k(\g,f)\on{-Mod},\quad
M\mapsto H^0_{DS,f}(M)\label{eq:DS-functor}
\end{align}
is not exact for a general nilpotent element $f$,
its restriction to $\on{KL}_k$ is exact:
\begin{Th}[\cite{Ara09b}]\label{Th;exact-kl}
Let $f$, $k$ be arbitrary.
We have
$H^{i}_{DS,f_{min}}(M)=0$
for any
 $M\in \on{KL}_k$ and $i\ne 0$.
Therefore,
the functor
$$\on{KL}_k\ra \W^k(\g,f_{min})\on{-Mod},\quad
M\mapsto H^0_{DS,f_{min}}(M)$$
is exact.
\end{Th}
In the case that $f$ is a principal nilpotent element,
Theorem \ref{Th;exact-kl} has been proved in \cite{FreGai07}
using Theorem \ref{Th:principal} below.

Unfortunately, the restriction of the quantized Drinfeld-Sokolv reduction functor
to
$\on{KL}_k$  does not produce all the irreducible highest weight representations of $\W^k(\g,f)$.
However, one can modify the functor \eqref{eq:DS-functor}
to the ^^ ^^ $-$"-reduction functor $H^0_{-,f}(?)$ (defined in \cite{FKW92})
to obtain the following result for the principal $W$-algebras.

\begin{Th}[\cite{Ara07}]\label{Th:principal}
 Let $f$ be a principal nilpotent element,
and 
let $k$ be an arbitrary complex number.
\begin{enumerate}
\item  
We have
$H^{i}_{-,f_{prin}}(M)=0$
for any
 $M\in \mc{O}_k$ and $i\in \Z\backslash \{0\}$.
Therefore,
the functor
$$\mc{O}_k\ra \W^k(\g,f_{prin})\on{-Mod},\quad
M\mapsto H^{i}_{-,f_{prin}}(M)$$
is exact.
\item 
For a weight $\lam$ of $\affg$ of level $k$, 
$H^{0}_{-,f_{prin}}(L(\lam))$
is zero or isomorphic to an
irreducible highest weight representation of
$\W^k(\g,f_{prin})$.
Moreover, any 
irreducible highest weight representation
the principal $W$-algebra $\W^k(\g,f_{prin})$ arises in this way.
 \end{enumerate}
\end{Th}

The vanishing and the irreducibility
statement of  Theorem \ref{Th:minimal} and Theorem \ref{Th:principal}
was conjectured for admissible representations of $\affg$ (see below)
in \cite{KacRoaWak03} and \cite{FKW92}, respectively.

In type $A$ we
 can derive the similar result as Theorem \ref{Th:principal}
for any nilpotent element $f$
by upgrading  the work of Brundan and Kleshchev \cite{BruKle08}
on the representation theory of finite W-algebras
to the affine setting
(\cite{Ara08-a}).
In particular,
the character of all  ordinary irreducible representations
of $\W^k(\mf{sl}_n,f)$ has been determined for a non-critical $k$.
\section{BRST reduction of associated varieties}
 Let $\W_k(\g,f)$ be the unique simple 
graded quotient of $\W^k(\g,f)$.
The associated variety  $X_{\W_k(\g,f)}$
is a subvariety of $X_{\W^k(\g,f)}=\Slo_f$,
which is invariant under the 
natural $\C^*$-action on $\Slo_f$ that contracts to the point $f\in \Slo_f$.
Therefore $\W_k(\g,f)$ is lisse if and only if $X_{\W_k(\g,f)}=\{f\}$.

By Theorem \ref{Th;exact-kl},
$\W_k(\g,f)$
is a quotient 
of the vertex algebra $H^0_{DS,f}(L_k(\g))$,
provided that it is nonzero.
 \begin{Th}[\cite{Ara09b}] \label{Th:IMRN}
For any $f\in \g$ and $k\in \C$ we have 
$$X_{H_{DS,f}^0(L_k(\g))}\cong X_{L_k(\g)}\cap \mathcal{S}_f.$$
Therefore,
 \begin{enumerate}
\item $H_{DS,f}^0(L_k(\g))\ne 0$ if and only if  $ \overline{G.f}\subset X_{L_k(\g)}$;
\item If  $X_{L_k(\g)}=\overline{G.f}$
then
$X_{H_{DS,f}^0(L_k(\g))}=\{f\}$.
Hence
$H_{DS,f}^0(L_k(\g))$ is lisse, and thus, so is its quotient
$\W_k(\g,f)$.
 \end{enumerate}
 \end{Th}
 Theorem \ref{Th:IMRN}
 can be regarded as a vertex algebra
 analogue
 of the corresponding result \cite{Los11,Gin08}
 for finite $W$-algebras.  
 
Note that 
if $L_k(\g)$ is integrable
we have
$H_{DS,f}^0(L_k(\g))=0$
  by \eqref{eq:integrability}.
Therefore we need to study more general representations of $\affg$
to obtain lisse $W$-algebras using Theorem \ref{Th:IMRN}.

Recall that the irreducible highest weight representation
$L(\lam)$ of $\affg$ is called {\em admissible} \cite{KacWak89}
(1) if $\lam$ is regular dominant, that is,
$\langle \lam+\rho,\alpha^{\vee}\rangle\not\in -\Z_{\geq 0}$ for any  $\alpha\in \Delta^{re}_+$,
and (2) $\mathbb{Q}\Delta(\lam)=\mathbb{Q}\Delta^{re}$.
Here $\Delta^{re}$  is the set of real roots of $\affg$,
$\Delta^{re}_+$ the set of positive real roots of $\affg$,
and $\Delta(\lam)=\{\alpha\in\Delta^{re}\mid \langle \lam+\rho,\alpha^{\vee}\rangle\in \Z\}$,
the set  of integral roots if $\lam$.
Admissible representations are (conjecturally all) modular invariant representations
of $\affg$,
that is,
 the characters of admissible representations are invariant under the natural action of $SL_2(\Z)$  (\cite{KacWak88}).
The simple affine vertex algebra $L_k(\g)$ is admissible as a $\affg$-module if and only if 
\begin{align}
k+h^{\vee}=\frac{p}{q},\quad p,q\in \mathbb{N},\ (p,q)=1,\ 
p\geq \begin{cases}
h^{\vee}&\text{if }(q,r^{\vee})=1,\\
h&\text{if }(q,r^{\vee})=r^{\vee}
\end{cases}
\label{eq:ad-number}
\end{align}
(\cite{KacWak08}).
Here 
$h$ is the Coxeter number of $\g$ and 
$r^{\vee}$ is the lacity of $\g$.
If this is the case
$k$ is called an {\em admissible number} for $\affg$ and  $L_k(\g)$ is called an {\em admissible affine vertex algebra}.

\begin{Th}[\cite{Ara09b}]\label{Th:ass-variety-of admissible}
Let  $L_k(\g)$ be an  admissible affine vertex algebra.
 \begin{enumerate}
\item (Feigin-Frenkel conjectrue)
We have $X_{L_k(\g)}\subset \mc{N}$. 
\item 
The variety $X_{L_k(\g)}$ is irreducible.
That is, there exists a nilpotent orbit $\mathbb{O}_k$ of $\g$
such that
$$X_{L_k(\g)}=\overline{\mathbb{O}_k}.$$
\item More precisely,
let $k$ be an admissible number of the form \eqref{eq:ad-number}.
Then 
\begin{align*}
X_{L_k(\g)}=\begin{cases}
\{x\in \g\mid (\on{ad}x)^{2q}=0\}&\text{if }(q,r^{\vee})=1,\\
\{x\in\g\mid \pi_{\theta_s}(x)^{2q/r^{\vee}}=0\}&\text{if }(q,r^{\vee})=r^{\vee},
\end{cases}
\end{align*}
where 
$\theta_s$ is the highest short root of $\g$
and
$\pi_{\theta_s}$ is the irreducible finite-dimensional representation 
of  $\g$ with highest weight $\theta_s$.
\end{enumerate}

\end{Th}

From Theorem \ref{Th:IMRN}
and Theorem \ref{Th:ass-variety-of admissible}
we immediately obtain the following assertion,
which  was (essentially) conjectured by Kac and Wakimoto \cite{KacWak08}.

\begin{Th}[\cite{Ara09b}]
\label{Th:lisse}
Let 
$L_k(\g)$ be an admissible affine vertex algebra,
and let
$f\in \mathbb{O}_k$.
Then the simple affine $W$-algebra
$\W_k(\g,f)$ is lisse.
\end{Th}
In the case that 
$X_{L_k(\g)}=\overline{G.f_{prin}}$,
the lisse $W$-algebras
obtained in Theorem \ref{Th:lisse}
is the {\em minimal series} principal $W$-algebras
 studied in \cite{FKW92}.
 In the case that $\g=\mf{sl}_2$,
these 
are exactly the minimal series Virasoro vertex algebras (\cite{FeuFuc84,BeiFeiMaz,Wan93}).
The Frenkel-Kac-Wakimoto conjecture 
states that these  minimal series principal $W$-algebras are rational\footnote{More generally,
all the lisse $W$-algebras $ \W_k(\g,f)$ that appear in Theorem \ref{Th:lisse}
are conjectured to be rational (\cite{KacWak08, Ara09b}).}.

\section{The rationality of minimal series principal W-algerbas}
An admissible affine vertex algebra
$L_k(\g)$ is called {\em non-degenearte} (\cite{FKW92})
if
$$X_{L_k(\g)}=\mc{N}=\overline{G.f_{prin}}.$$
If this is the case $k$ is called a {\em non-degenerate admissible number} for $\affg$.
By theorem \ref{Th:ass-variety-of admissible} (iii),
^^ ^^ most"  admissible affine vertex algebras are non-degenerate.
More precisely,
an admissible number $k$
of the form \eqref{eq:ad-number}
is non-degenerate if and only if $q\geq \begin{cases} h&\text{if }(q,r^{\vee})=1,\\ r^{\vee}
{}^Lh^{\vee}&\text{if }(q,r^{\vee})=r^{\vee}\end{cases}$,
where ${}^Lh^{\vee}$ is the dual Coxeter number of the Langlands dual Lie algebra ${}^L\g$.
For a non-degenerate admissible number $k$,
the simple principal $W$-algebra $\W_k(\g,f_{prin})$ is lisse by Theorem \ref{Th:lisse}.

The following assertion settles 
the  Frenkel-Kac-Waimoto conjecture \cite{FKW92}
in full generality.
 \begin{Th}[\cite{A2012Dec}]\label{Th:minimal-series}
Let $k$ be a non-degenerate admissible number.
Then
  the simple principal $W$-algebra  $\W_k(\g,f_{prin})$ is rational.
\end{Th}

The proof of Theorem \ref{Th:minimal-series} based on
Theorem \ref{Th:principal}, Theorem \ref{Th:lisse},
and the following assertion on admissible affine vertex algebras,
which  was conjectured by  Adamovi\'{c} and Milas \cite{AdaMil95}.
\begin{Th}[\cite{A12-2}]
Let 
 $L_k(\g)$ be an admissible affine vertex algebra.
 Then $L_k(\g)$ is rational in the category $\mc{O}$,
 that is,
 any $L_k(\g)$-module that belongs to  $\mc{O}$ is completely reducible.
\end{Th}

The following assertion,
that   has been widely believed 
since \cite{KacWak90},
gives a yet another realization of minimal series principal $W$-algebras.
\begin{Th}[\cite{ACL17}]
\label{Th:GKO}
Let $\g$ be simply laced.
For  an admissible affine vertex algebra
$L_k(\g)$,
 the  vertex algebra 
$ (L_k(\g)\otimes L_1(\g))^{\g[t]}$ is isomorphic to a
minimal series principal $W$-algebra.
Conversely,  any minimal series principal $W$-algebra  associated with $\g$  appears in this way.
\end{Th}
In the case that $\g=\mf{sl}_2$ and
$k$ is a non-negative integer, 
the statement of Theorem \ref{Th:GKO}
is well-known as the GKO construction of the discrete series 
of the Virasoro vertex algebras \cite{GodKenOli86}.
Some partial results have been obtained previously in \cite{ALYcoset, AJiang}.
From Theorem \ref{Th:GKO}, it follows that 
the minimal series principal $W$-algebra $\W_{p/q-h^{\vee}}(\g,f_{prin})$  of  ADE type is 
unitary, that is, any simple $\W_{p/q-h^{\vee}}(\g,f_{prin})$-module is 
unitary in the sense of \cite{DonLin14},
if and only if $|p-q|=1$.
\section{Four-dimensional $N=2$ superconformal algebras,
 Higgs branch conjecture
 and the class $\mathcal{S}$ chiral algebras}
 In the study of four-dimensional $N=2$ superconformal field theories in physics,
Beem,  Lemos, Liendo, Peelaers, Rastelli,
and van Rees  \cite{BeeLemLie15} have constructed a remarkable 
map
\begin{align}
\Phi: \{\text{4d $N=2$ SCFTs}\}\ra \{\text{vertex algebras}\}
\label{eq:the-map-phi}
\end{align}
such that,
among other things,
the character of the vertex algebra
$\Phi(\mc{T})$ coincides  with the {\em Schur index}
of the corresponding   4d $N=2$ SCFT $\mc{T}$,
which is an important invariant of the theory $\mc{T}$. 

How do vertex algebras coming from 4d $N=2$ SCFTs look like?
According to \cite{BeeLemLie15}, we have
$$c_{2d}=-12c_{4d},$$
where $c_{4d}$ and $c_{2d}$ are central charges of the 4d $N=2$ SCFT
and the corresponding vertex algebra, respectively.
Since the central charge is positive for a unitary theory,
this implies that the vertex algebras obtained by $\Phi$ are never unitary.
In particular integrable affine vertex algebras never appear by this correspondence.

The main examples of vertex algebras
considered in  \cite{BeeLemLie15}
are the simple affine vertex algebras $L_k(\g)$ of types 
$D_4$, $F_4$, $E_6$, $E_7$, $E_8$ at level $k=-h^{\vee}/6-1$,
which are non-rational,
non-admissible affine vertex algebras 
studied in \cite{AM15} by a different motivation.
One can find more examples in the literature, see e.g.
 \cite{Beem:2015yu,BN1,
CorSha16,
BN2,
Dan,SXY,BLN}.

Now, there is another important invariant of a 4d $N=2$ SCFT $\mc{T}$,
called the {\em Higgs branch},
which we denote by $Higgs_{\mc{T}}$.
The Higgs branch $Higgs_{\mc{T}}$ 
is an affine  algebraic variety
that has the hyperK\"{a}hler structure in its smooth part.
In particular,
$Higgs_{\mc{T}}$  is a (possibly singular)
 symplectic  variety.
 
 Let $\mc{T}$ be one of the  4d $N=2$ SCFTs studied in  \cite{BeeLemLie15}
such that  that $\Phi(\mc{T})=L_k(\g)$ with 
$k=h^{\vee}/6-1$
for types 
$D_4$, $F_4$, $E_6$, $E_7$, $E_8$ as above. 
It is known that $Higgs_{\mc{T}}=\overline{\mathbb{O}_{min}}$,
and this equals  \cite{AM15} to the
the associated variety $X_{\Phi(\mc{T})}$.
It is expected that this is not just a coincidence.

\begin{Conj}[Beem and Rastelli \cite{BeeRas}]\label{Conj:Beem and Rastell}
For a 4d $N=2$ SCFT $\mc{T}$, we have
\begin{align*}
Higgs_{\mc{T}}=X_{\Phi(\mc{T})}.
\end{align*}
\end{Conj}
So we are expected to  recover the Higgs branch of a 4d $N=2$ SCFT
from the corresponding vertex algebra, which is 
a purely algebraic object.

We note that Conjecture \ref{Conj:Beem and Rastell} is a physical conjecture
since the Higgs branch is not  a  mathematically defined object at the moment.
 The Schur index is not a mathematically defined object  either.
However, in view of \eqref{eq:the-map-phi}  and Conjecture \ref{Conj:Beem and Rastell},
one can try to define both Higgs branches and Schur indices
of 4d $N=2$ SCFTs
using vertex algebras.
We note that there is a close relationship between the Higgs branches of 4d $N=2$ SCFTs
and the {\em Coulomb branches} of three-dimensional  $N=4$ gauge theories whose mathematical definition has been given by 
Braverman, Finkelberg and Nakajima \cite{BFN16}.

Although Higgs branches are symplectic varieties,
the associated variety  $X_V$  of a vertex algebra  $V$ is only a Poisson variety in general.
A vertex algebra $V$ is called {\em quasi-lisse} (\cite{Arakawam:kq}) if $X_V$ has only finitely many symplectic leaves.
If this is the case symplectic leaves in $X_V$  are algebraic (\cite{Brown-Gordon}).
Clearly,
lisse vertex algebras are quasi-lisse.
The simple affine vertex algebra $L_k(\g)$
is quasi-lisse if and only if $X_{L_k(\g)}\subset \mc{N}$.
In particular, admissible affine vertex algebras are quasi-lisse.
See \cite{AM15, AMsheets, AMirr} for more examples of quasi-lisse vertex algebras.
Physical intuition expects that
vertex algebras that come from 4d $N=2$ SCFTs via the map $\Phi$   are quasi-lisse.

By extending Zhu's argument \cite{Zhu96}
using a theorem of Etingof and Schelder  \cite{EtiSch10},
we obtain the following assertion.
\begin{Th}[\cite{Arakawam:kq}]\label{Th:AK}
Let $V$ be a quasi-lisse $\Z_{\geq 0}$-graded conformal vertex algebra such that $V_0=\C$.
Then there only finitely many simple
 ordinary   $V$-modules.
 Moreover, 
for  a finitely generated ordinary $V$-module
$M$,
the character $\chi_M(q)$  satisfies a {\em modular linear differential equation}.
\end{Th}
Since the space of solutions of a modular linear differential equation is invariant under the action of $SL_2(\Z)$,
 Theorem \ref{Th:AK} implies that a quasi-lisse vertex algebra possesses a certain  modular invariance property,
although we do not claim that the normalized characters of ordinary $V$-modules span the space of the solutions.
An important consequence of  Theorem \ref{Th:AK}  
is that   the Schur indices of 4d $N=2$ SCFTs
have some modular invariance property.
This is something that has been conjectured by physicists (\cite{BeeRas}).

\smallskip

There is a distinct class of 
four-dimensional $N=2$ superconformal field theories
called the {\em theory of class $\mc{S}$} (\cite{Gai12,GaiMooNei13}),
where $\mc{S}$ stands for $6$.
The vertex algebras obtained from the theory of class $\mc{S}$
is called the {\em chiral algebras of class $\mc{S}$} (\cite{Beem:2015yu}).
The Moore-Tachikawa conjecture \cite{MooTac12},
which was recently proved in \cite{BraFinNak17},
describes the Higgs branches of the  theory of class $\mc{S}$
in terms of two-dimensional topological quantum field theories
mathematically.

Let $\mathbb{V}$ be the category of vertex algebras,
whose objects are complex semisimple groups,
and 
$\on{Hom}(G_1,G_2)$ is the isomorphism classes
of conformal vertex algebras $V$ with a vertex algebra homomorphism
\begin{align*}
V^{-h_1^{\vee}}(\g_1)\otimes V^{-h_2^{\vee}}(\g_2)\ra V\label{eq:quatized-momentmap}
\end{align*}
such that the action of $\g_1[t]\+ \g_2[t]$ on $V$ is locally finite.
Here $\g_i=\on{Lie}(G_i)$ and $h_i^{\vee}$ is the dual Coxeter number of $\g_i$
in the case that $\g_i$ is simple.
If $\g_i$  is not simple
we understand $V^{-h_i^{\vee}}(\g_i)$ to be the tensor product of 
the
 critical level universal affine vertex algebras corresponding to  all simple components of $\g_i$.
 The composition $V_1\circ V_2$ of $V_1\in \on{Hom}(G_1,G_2)$ 
 and $V_2\in \on{Hom}(G_1,G_2)$ 
 is given by the relative semi-infinite cohomology
 \begin{align*}
V_1\circ V_2=H^{\frac{\infty}{2}+0}(\affg_2, \g_2,V_1\otimes V_2),
\end{align*}
where $\affg_2$ denotes the direct sum of the affine Kac-Moody algebra
associated with the simple components of $\g_2$.
By a result of \cite{ArkGai02},
one finds that
the identity morphism $\on{id}_G$ is the algebra $\mc{D}_G^{ch}$ of {\em chiral differential operators}
on $G$ (\cite{MalSchVai99,BeiDri04}) at the critical level,
whose associated variety is canonically isomorphic to $T^*G$.

The following theorem,
which was conjectured  in  \cite{BeeLemLie15}
(see  \cite{TachikawaVOA,TachikawaICM}) for  mathematical expositions),
describes the chiral algebras of class $\mc{S}$.
\begin{Th}[\cite{Achiral}]\label{Th:chiral-algebras-of-class-S}
Let $\mathbb{B}_2$ the category of $2$-bordisms.
There exists a unique monoidal functor
$\eta_G:\mathbb{B}_2\ra \mathbb{V}$
which sends 
(1) the object $S^1$ to $G$, 
(2) the cylinder,
which is the identity morphism $\on{id}_{S^1}$,
to the 
 identity morphism $\on{id}_G=\mc{D}_G^{ch}$,
 and (3) the cap to $H^0_{DS,f_{prin}}(\mc{D}_G^{ch})$.
 Moreover,
we have  $X_{\eta_G(B)}\cong \eta^{BFN}_G(B)$   for any 2-bordism $B$,
 where $\eta^{BFN}_G$ is the functor form $\mathbb{B}_2$  to the category of
 symplectic variaties constructed in \cite{BraFinNak17}.
\end{Th}
In view of \cite{MooTac12,BeeLemLie15},
the last assertion of Theorem \ref{Th:chiral-algebras-of-class-S}
confirms 
the Higgs branch conjecture (Conjecture \ref{Conj:Beem and Rastell}) for the theory of class $\mc{S}$.

\newcommand{\etalchar}[1]{$^{#1}$}

\end{document}